\documentclass[a4paper,11pt]{amsart}

\usepackage{graphicx}
\usepackage{mathptmx}
\usepackage{amsmath}
\usepackage{amssymb}
\usepackage{enumitem}
\usepackage{xcolor}

\newmuskip\pFqmuskip

\newcommand*\pFq[6][8]{%
  \begingroup 
  \pFqmuskip=#1mu\relax
  \mathcode`=\string"8000
  \begingroup\lccode`\~=`\,
  \lowercase{\endgroup\let~}\pFqcomma
  F^{#2}_{#3}{\left(\genfrac..{0pt}{}{#4}{#5}\bigg|#6\right)}%
  \endgroup
}
\newcommand{\pFqcomma}{\mskip\pFqmuskip}

\newcommand{\ap}{a^{\dagger}}
\newcommand{\ba}{aa^{\dagger}}
\newcommand{\fa}{a^{\dagger}a}

\newtheorem{theorem}{Theorem}[section]

\newtheorem{proposition}[theorem]{Proposition}

\begin{document}

\title[]{Recurrence relations for degenerate Bell and Dowling polynomials via Boson operators}

\author{Taekyun  Kim}
\address{Department of Mathematics, Kwangwoon University, Seoul 139-701, Republic of Korea}
\email{tkkim@kw.ac.kr}

\author{Dae San  Kim}
\address{Department of Mathematics, Sogang University, Seoul 121-742, Republic of Korea}
\email{dskim@sogang.ac.kr}

\subjclass[2010]{11B73; 11B83}
\keywords{degenerate Bell polynomials; degenerate Dowling polynomials; degenerate $r$-Dowling polynomials; Boson operator}

\begin{abstract}
Spivey found a recurrence relation for the Bell numbers by using combinatorial method.
The aim of this paper is to derive Spivey's type recurrence relations for the degenerate Bell polynomials and the degenerate Dowling polynomials by using the boson annihilation and creation operators satisfying the commutation relation $aa^{+}-a^{+}a=1$. In addition, we derive a Spivey's type recurrence relation for the $r$-Dowling polynomials.
\end{abstract}

\maketitle

\markboth{\centerline{\scriptsize Recurrence relations for degenerate Bell and Dowling polynomials via Boson operators}}
{\centerline{\scriptsize  Taekyun Kim and Dae San Kim}}

\section{Introduction}
In [32], Spivey discovered by using combinatorial method a recurrence relation for the Bell numbers
\begin{equation}
\phi_{n+m}=\sum_{k=0}^{n}\sum_{j=0}^{m}{m \brace j}\binom{n}{k}j^{n-k}\phi_{k},\quad (\mathrm{see}\ [11,21,32]), \label{-1}
\end{equation}
which generalizes both of the following well-known expressions for those numbers (see \eqref{4}, \eqref{5}):
\begin{equation*}
\phi_{n}=\sum_{j=0}^{n}{n \brace j},\quad \phi_{n+1}=\sum_{k=0}^{n}\binom{n}{k}\phi_{k}.
\end{equation*} \par
Degenerate versions of many special polynomials and numbers have been explored by using various methods, including quantum mechanics, combinatorial methods, generating functions, umbral calculus, $p$-adic analysis, probability theory, differential equations, operator theory and analytic number theory. \par
The aim of this paper is to derive Spivey's type recurrence relations for the degenerate Bell polynomials, the degenerate Dowling polynomials and the degenerate $r$-Dowling polynomials by using the boson annihilation and creation operators satisfying the commutation relation $aa^{+}-a^{+}a=1$. \par
In more detail, the outline of this paper is as follows.
In Section 2, we derive a recurrence relation for degenerate Bell polynomials
\begin{equation}
\phi_{n+m,\lambda}(x)=\sum_{j=0}^{m}\sum_{k=0}^{n}{m \brace j}_{\lambda}\binom{n}{k}(j-m\lambda)_{n-k,\lambda}x^{j}\phi_{k,\lambda}(x) \label{0}
\end{equation}
from the normal ordering of `degenerate $n$-th power' of the number operator $\hat{n}=\ap a$
\begin{equation*}
\big(\ap a\big)_{n,\lambda}=\sum_{k=0}^{n}{n \brace k}_{\lambda}\big(\ap\big)^{k}a^{k},
\end{equation*}
in terms of boson annihilation and creation operators $a$ and $\ap$ satisfying the commutation relation with $a\ap-\ap a=1$ (see \eqref{8}, \eqref{10}).
Note that taking $x=1$ and letting $\lambda \rightarrow 0$ in \eqref{0} yield the Spivey's recurrence relation for the Bell numbers in \eqref{-1}. \par
In Section 3, we deduce recurrence relations for the degenerate Dowling polynomials $D_{m,\lambda}(n,x)$ (see Theorem 3.1) and the degenerate $r$-Dowling polynomials $D_{m,\lambda}^{(r)}(n,x)$ (see Theorem 3.2), respectively from the normal orderings of degenerate $n$-th powers of the operators $m\ap a+1$  (see \eqref{34}) and $m\ap a+r$ (see \eqref{35}). For the rest of this section, we recall the facts that are needed  throughout this paper. \par

\vspace{0.1in}

For any nonzero $\lambda\in\mathbb{R}$, the degenerate exponentials are defined by
\begin{equation}
e_{\lambda}^{x}(t)=(1+\lambda t)^{\frac{x}{\lambda}}=\sum_{k=0}^{\infty}(x)_{k,\lambda}\frac{t^{k}}{k!},\quad e_{\lambda}(t)=e_{\lambda}^{1}(t), \label{1}
\end{equation}
where
\begin{equation}
(x)_{0,\lambda}=1,\quad (x)_{n,\lambda}=x(x-\lambda)(x-2\lambda)\cdots \big(x-(n-1)\lambda\big),\quad (n\ge 1),\quad (\mathrm{see}\ [13-27]). \label{2}
\end{equation}
The degenerate Stirling numbers of the second kind are given by
\begin{equation}
(x)_{n,\lambda}=\sum_{k=0}^{n}{n \brace k}_{\lambda}(x)_{k},\quad (n\ge 0),\quad (\mathrm{see}\ [9,13-27,30]).\label{3}
\end{equation}
Note that $\displaystyle\lim_{\lambda\rightarrow 0}{n \brace k}_{\lambda}={n \brace k}\displaystyle$ are the Stirling numbers of the second kind defined by
\begin{equation}
x^{n}=\sum_{k=0}^{n}{n \brace k}(x)_{k},\quad (n\ge 0),\quad (\mathrm{see}\ [1-33]). \label{4}
\end{equation} \par
The Bell polynomials are given by
\begin{equation}
\sum_{n=0}^{\infty}\phi_{n}(x)\frac{t^{n}}{n!}=e^{x(e^{t}-1)},\quad (\mathrm{see}\ [8,10,12,26,31]). \label{5}
\end{equation}
When $x=1$, $\phi_{n}=\phi_{n}(1),\ (n\ge 0)$, are called the Bell numbers.
The degenerate Bell polynomials are defined by
\begin{equation}
\phi_{n,\lambda}(x)=\sum_{k=0}^{n}{n \brace k}_{\lambda}x^{k},\quad (n,k\ge 0),\quad (\mathrm{see}\ [16-19]). \label{7}
\end{equation}
From \eqref{7}, we note that (see \eqref{1})
\begin{equation*}
e^{x(e_{\lambda}(t)-1)}=\sum_{n=0}^{\infty}\phi_{n,\lambda}(x)\frac{t^{n}}{n!},\quad (\mathrm{see}\ [14-18,23,24]).
\end{equation*}
Note that
\begin{displaymath}
\lim_{\lambda\rightarrow 0}\phi_{n,\lambda}(x)=\phi_{n}(x),\quad (n\ge 0).
\end{displaymath}
When $x=1,\ \phi_{n,\lambda}=\phi_{n,\lambda}(1)$ are called the degenerate Bell numbers. \par
We recall that the boson annihilation and creation operators $a$ and $\ap$ satisfy the commutation relation given by
\begin{equation}
[a,\ap]=\ba-\fa=1,\quad (\mathrm{see}\ [4,5,11-13,16-18,23-25,27,29]).\label{8}	
\end{equation}
In the corresponding Fock space, we note that
\begin{equation}
\ap|k\rangle=\sqrt{k+1}|k+1\rangle,\quad a|k\rangle=\sqrt{k}|k-1\rangle,\quad (\mathrm{see}\ [11-13,16-18,23-25,27,29]).\label{9}
\end{equation}
By \eqref{9}, we get $\ap a|k\rangle=k|k\rangle$. The coherent states $|z\rangle$, where $z$ is complex number, satisfy $a|z\rangle=z|z\rangle,\ \langle z|z\rangle=1$ and $\langle z|a=\langle z|\overline{z}$. We recall that the harmonic oscillator has Hamiltonian $\hat{n}=\ap a$ (neglecting the zero point energy) and the usual eigenstates $|n\rangle$ (for $n\in\mathbb{N}$) satisfying $\hat{n}|n\rangle=n|n\rangle$ and $\langle m|n\rangle=\delta_{m,n}$, where $\delta_{m,n}$ is the Kronecker's symbol (see [16]). It is known that the standard basis commutation relation $a\ap-\ap a=1$ can be considered formally, in a suitable space of functions, by letting $a=\frac{d}{dx}$ and $\ap=x$. \par
We note that (see \eqref{4}, \eqref{15})
\begin{equation*}
\big(\ap a\big)^{k}=\sum_{l=0}^{k}{k \brace l}\big(\ap)^{l}a^{l},\quad (k\ge 0),\quad (\mathrm{see}\ [11-13,16-18,23-25,27]).	
\end{equation*}
The normal ordering of degenerate $n$-th power of the number operator $\hat{n}=\ap a$ is given in terms of boson operators $a$ and $\ap$ by (see \eqref{2}, \eqref{3}, \eqref{15})
\begin{equation}
\big(\ap a\big)_{n,\lambda}=\sum_{k=0}^{n}{n \brace k}_{\lambda}\big(\ap\big)^{k}a^{k},\label{10}
\end{equation}
and
\begin{displaymath}
\big(\ap\big)^{n}a^{n}=\sum_{k=0}^{n}S_{1,\lambda}(n,k)\big(\ap a\big)_{k,\lambda},\quad (n\ge 0),\ (\mathrm{see}\ [16]),
\end{displaymath}
where $S_{1,\lambda}(n,k)$ are the degenerate Stirling numbers of the first kind defined by
\begin{equation}
(x)_{n}=\sum_{k=0}^{n}S_{1,\lambda}(n,k)(x)_{k,\lambda},\quad (n\ge 0),\quad (\mathrm{see}\ [16,17,19,22]).\nonumber
\end{equation}

\section{Recurrence relation for degenerate Bell polynomials via Boson operator}
In this section, we derive a Spivey's type recurrence relation for the degenerate Bell polynomials via boson operator. \par
From \eqref{8}, we note that
\begin{align}
\big[a,\hat{n}\big]&=a\hat{n}-\hat{n}a=a\big(\ap a)-(\ap a)a\label{11}	\\
&=\big(a\ap-\ap a\big)a=a\nonumber
\end{align}
and
\begin{align}
\big[\hat{n},a^{+}\big]&=\hat{n}\ap-\ap \hat{n}=\big(\ap a\big)\ap-\ap\big(\ap a\big) \label{12} \\
&=\ap\big(a\ap-\ap a\big)=\ap. \nonumber
\end{align}
For $k\in\mathbb{N}$, by \eqref{8}, we get
\begin{align}
\big[a,\big(\ap\big)^{k}\big]&=a\big(\ap\big)^{k}-\big(\ap\big)^{k}a \label{13}\\
&=\big(a\big(\ap\big)^{k-1}-\big(\ap\big)^{k-1}a\big)\ap+\big(\ap)^{k-1}\big(a\ap-\ap a\big)\nonumber\\
&=\big[a,\big(\ap\big)^{k-1}\big]\ap+\big(\ap\big)^{k-1}\big[a,\ap\big]\nonumber\\
&=\big[a,\big(\ap\big)^{k-1}\big]\ap+\big(\ap\big)^{k-1}=\big[a,\big(\ap\big)^{k-2}\big]\big(\ap\big)^{2}+2\big(\ap\big)^{k-1} \nonumber\\
&=\cdots\nonumber\\
&=\big[a,\ap\big]\big(\ap\big)^{k-1}+(k-1)\big(\ap\big)^{k-1}=k\big(\ap\big)^{k-1}.\nonumber
\end{align}
In a similar manner, we can show that
\begin{equation}
\big[a^k,\ap\big]=ka^{k-1}, \label {14}
\end{equation}
which in turn implies that
\begin{equation}
(\ap a)_{k}=(\ap)^{k}a^{k}. \label{15}
\end{equation} \par
For $z\in\mathbb{C}$, the coherent states $|z\rangle$ satisfy $a|z\rangle=z|z\rangle$. Note that $a|0\rangle=0|0\rangle=0$. \\
Thus, we have
\begin{align}
\big[a,\big(\ap\big)^{k}\big]\big|0\big\rangle &=a\big(\ap)^{k}-\big(\ap\big)^{k}a\big|0\big\rangle=a\big(\ap\big)^{k}\big|0\big\rangle-\big(\ap\big)^{k}a\big|0\big\rangle \label{16}\\
&=a\big(a^{+}\big)^{k}\big|0\big\rangle,\quad (k\in\mathbb{N}).\nonumber	
\end{align}
From \eqref{13} and \eqref{16}, we note that
\begin{equation}
a\big(\ap\big)^{k}\big|0\big\rangle=\big[a,\big(\ap\big)^{k}\big]\big|0\big\rangle=k\big(\ap\big)^{k-1}\big|0\big\rangle.\label{17}
\end{equation}
Therefore, by \eqref{11}, \eqref{12}, \eqref{13}, \eqref{14} and \eqref{17}, we obtain the following proposition.
\begin{proposition}
For $k\in\mathbb{N}$, we have
\begin{align*}
&\big[a,\big(\ap\big)^{k}\big]=k\big(\ap\big)^{k-1},\quad \big[a^k,\ap\big]=ka^{k-1}, \\
&a\big(\ap\big)^{k}\big|0\big\rangle=\big[a,\big(\ap\big)^{k}\big]\big|0\big\rangle=k\big(\ap\big)^{k-1}\big|0\big\rangle.
\end{align*}
In addition,
\begin{displaymath}
\big[a,\hat{n}\big]=a\quad \mathrm{and}\quad [\hat{n},\ap]=\ap.
\end{displaymath}
\end{proposition}
From \eqref{10}, we have
\begin{equation}
\big(\ap a\big)_{n,\lambda}e^{\ap}\big|0\big\rangle=\sum_{k=0}^{n}{n \brace k}_{\lambda}\big(\ap\big)^{k}a^{k}e^{\ap}\big|0\big\rangle.\label{18}	
\end{equation}
Now, we observe from \eqref{17} that
\begin{align}
a^{k}e^{\ap}\big|0\big\rangle &=\sum_{l=0}^{\infty}\frac{1}{l!}a^{k}\big(\ap\big)^{l}\big|0\big\rangle=\sum_{l=1}^{\infty}\frac{1}{(l-1)!}a^{k-1}\big(\ap\big)^{l-1}\big|0\big\rangle\label{19} \\
&=\sum_{l=0}^{\infty}\frac{1}{l!}a^{k-1}\big(\ap\big)^{l}\big|0\big\rangle=\sum_{l=1}^{\infty}\frac{1}{(l-1)!}a^{k-2}\big(\ap\big)^{l-1}\big|0\big\rangle\nonumber\\
&=\sum_{l=0}^{\infty}\frac{1}{l!}a^{k-2}\big(\ap\big)^{l}\big|0\big\rangle=\cdots\nonumber\\
&=\sum_{l=1}^{\infty}\frac{1}{(l-1)!}\big(\ap\big)^{l-1}\big|0\big\rangle=\sum_{l=0}^{\infty}\frac{1}{l!}\big(\ap\big)^{l}\big|0\big\rangle=e^{\ap}\big|0\big\rangle.\nonumber
\end{align}
By \eqref{3}, \eqref{7}, \eqref{18} and \eqref{19}, we get
\begin{align}
\big(\ap a\big)_{n,\lambda}e^{\ap}\big|0\big\rangle &=\sum_{k=0}^{n}{n\brace k}_{\lambda}\big(\ap\big)^{k}a^{k}e^{\ap}\big|0\big\rangle \label{20}\\
&=\sum_{k=0}^{n}{n \brace k}_{\lambda}\big(\ap\big)^{k}e^{\ap}\big|0\big\rangle=\phi_{n,\lambda}(\ap)e^{\ap}\big|0\big\rangle.\nonumber
\end{align}
Therefore, by \eqref{19} and \eqref{20}, we obtain the following theorem.
\begin{theorem}
For $n\ge 0$, we have
\begin{equation*}
a^{k}e^{\ap}\big|0\big\rangle=e^{\ap}\big|0\big\rangle, \quad   \big(\ap a\big)_{n,\lambda}e^{\ap}\big|0\big\rangle=\phi_{n,\lambda}(\ap)e^{\ap}|0\big\rangle.
\end{equation*}
\end{theorem}
From \eqref{13}, we note that
\begin{equation}
a\big(\ap\big)^{k}=\big(\ap\big)^{k}a+k\big(\ap\big)^{k-1},\quad (k\in\mathbb{N}). \label{21}
\end{equation}
By \eqref{21}, we get
\begin{align}
\hat{n}\big(\ap)^{k}&=\big(\ap a\big)(\ap)^{k}=\ap\big(a(\ap)^{k}\big)=\ap\Big(\big(\ap\big)^{k}a+k\big(\ap\big)^{k-1}\Big) \label{22}\\
&=\big(\ap\big)^{k+1}a+k\big(\ap\big)^{k}=\big(\ap\big)^{k}\big(\ap a+k\big).\nonumber	
\end{align}
Thus, by \eqref{22}, we get
\begin{align}
&\big(\ap a-m\lambda\big)_{n,\lambda}\big(\ap\big)^{j}=\big(\ap a-m\lambda\big)\big(\ap a-m\lambda-\lambda\big)\cdots \label{23}\\
&\quad\times \big(\ap a-m\lambda-(n-1)\lambda\big)\big(\ap\big)^{j}\nonumber\\
&=\big(\ap a-m\lambda\big)\big(\ap a-m\lambda-\lambda\big)\cdots\big(\ap a-m\lambda-(n-2)\lambda\big)\nonumber\\
&\qquad\times \big(\ap\big)^{j}\big(\ap a-m\lambda-(n-1)\lambda+j\big) \nonumber\\
&=\big(\ap a-m\lambda\big)\big(\ap a-m\lambda-\lambda\big) \cdots\big(\ap a-m\lambda-(n-2)\lambda\big)\nonumber\\
&\quad\times \big(a^{+}\big)^{j}\big(\ap a-m\lambda-(n-2)\lambda+j\big)\big(\ap a-m\lambda-(n-1)\lambda+j\big)\nonumber\\
&=\cdots\nonumber\\
&=\big(\ap\big)^{j}\big(\ap a-m\lambda+j\big)_{n,\lambda},\nonumber
\end{align}
where $n,j$ are nonnegative integers. \\
From \eqref{22} and \eqref{23}, the next theorem follows.
\begin{theorem}
Let $k$ be a nonnegative integer. For $m,n\ge 0$, we have
\begin{align*}
&\hat{n}\big(\ap\big)^{k}=\big(\ap a\big)\big(\ap\big)^{k}=\big(\ap\big)^{k}\big(\ap a+k\big), \\
&\big(\ap a-m\lambda\big)_{n,\lambda}\big(\ap\big)^{k}=\big(\ap\big)^{k}\big(\ap a-m\lambda+k\big)_{n,\lambda},\\
& m\big(\fa)\big(\ap\big)^{j}=m\big(\ap\big)^{j}(\ba+j)=\big(\ap\big)^{j}\big(m\fa+mj\big).
\end{align*}
\end{theorem}
From \eqref{10} and Theorem 2.3, we note that
\begin{align}
\big(\ap a\big)_{m+n,\lambda}&=\big(\ap a-m\lambda\big)_{n,\lambda}\big(\ap a\big)_{m,\lambda}\label{24}  \\
&=\big(\ap a-m\lambda\big)_{n,\lambda}\sum_{j=0}^{m}{m \brace j}_{\lambda}\big(\ap\big)^{j}a^{j} \nonumber \\
&=\sum_{j=0}^{m}{m \brace j}_{\lambda}\big(\ap a-m\lambda\big)_{n,\lambda}\big(\ap\big)^{j}a^{j}\nonumber \\
&=\sum_{j=0}^{m}{m \brace j}_{\lambda}\big(\ap\big)^{j}\big(\ap a+j-m\lambda\big)_{n,\lambda}a^{j}.\nonumber
\end{align}
Therefore, by \eqref{24}, we obtain the following theorem.
\begin{theorem}
For $m,n\ge 0$, we have
\begin{equation}
\big(\ap a\big)_{m+n,\lambda}=\sum_{j=0}^{m}{m \brace j}_{\lambda}\big(\ap\big)^{j}\big(\ap a+j-m\lambda\big)_{n,\lambda}a^{j}. \label{25}	
\end{equation}
\end{theorem}
We recall the following identity:
\begin{equation}
(x+y)_{n,\lambda}=\sum_{k=0}^{n}\binom{n}{k}(x)_{n-k,\lambda}(y)_{k,\lambda},\quad (n\ge 0). \label{26}
\end{equation}
From \eqref{25} and \eqref{26}, we note that
\begin{align}
\big(\fa\big)_{m+n,\lambda}&=\sum_{j=0}^{m}{m \brace j}_{\lambda}\big(\ap\big)^{j}\big(\ap a+j-m\lambda\big)_{n,\lambda}a^{j}\label{27}\\
&=\sum_{j=0}^{m}\sum_{k=0}^{n}{m \brace j}_{\lambda}\binom{n}{k}\big(j-m\lambda\big)_{n-k,\lambda}\big(\ap\big)^{j}\big(\ap a\big)_{k,\lambda}a^{j}.\nonumber
\end{align}
Thus, by Theorem 2.2 and \eqref{27}, we get
\begin{align}
&\phi_{n+m,\lambda}\big(\ap\big)e^{\ap}\big|0\big\rangle=\big(\ap a\big)_{n+m,\lambda}e^{\ap}\big|0\big\rangle 	\label{28}\\
&=\sum_{j=0}^{m}\sum_{k=0}^{n}{m \brace j}_{\lambda}\binom{n}{k}(j-m\lambda)_{n-k,\lambda}\big(\ap\big)^{j}\big(\ap a\big)_{k,\lambda}a^{j}e^{\ap}\big|0\big\rangle \nonumber\\
&=\sum_{j=0}^{m}\sum_{k=0}^{n}{m \brace j}_{\lambda}\binom{n}{k}(j-m\lambda)_{n-k,\lambda}\big(\ap\big)^{j}\big(a\ap\big)_{k,\lambda}e^{\ap}\big|0\big\rangle \nonumber\\
&=\sum_{j=0}^{m}\sum_{k=0}^{n}{m \brace j}_{\lambda}\binom{n}{k}(j-m\lambda)_{n-k,\lambda}\big(\ap\big)^{j}\phi_{k,\lambda}\big(\ap\big)e^{\ap}\big|0\big\rangle.\nonumber
\end{align}
Therefore, by \eqref{28}, we obtain the following theorem.
\begin{theorem}
For $m,n\ge 0$, we have
\begin{equation*}
\phi_{n+m,\lambda}\big(\ap\big)=\sum_{j=0}^{m}\sum_{k=0}^{n}{m \brace j}_{\lambda}\binom{n}{k}(j-m\lambda)_{n-k,\lambda}\big(\ap\big)^{j}\phi_{k,\lambda}\big(\ap\big),
\end{equation*} \par
and hence
\begin{equation*}
\phi_{n+m,\lambda}(x)=\sum_{j=0}^{m}\sum_{k=0}^{n}{m \brace j}_{\lambda}\binom{n}{k}(j-m\lambda)_{n-k,\lambda}x^{j}\phi_{k,\lambda}(x). 	
\end{equation*}
\end{theorem}

\section{Recurrence relation for degenerate Dowling polynomials via boson operator}

In this section, we derive Spivey's type recurrence relations for the degenerate Dowling and $r$-Dowling polynomials via boson operators.  \par
For $m\in\mathbb{N}$, the degenerate Whitney numbers of the second kind are defined by
\begin{equation}
(mx+1)_{n,\lambda}=\sum_{k=0}^{n}W_{m,\lambda}(n,k)m^{k}(x)_{k},\quad (n\ge 0),\quad (\mathrm{see}\ [17]).\label{30}
\end{equation}
Note that
\begin{displaymath}
W_{1,\lambda}(n,k)={n+1 \brace k+1}_{\lambda}+\lambda n{n \brace k+1}_{\lambda},\quad (n\ge k\ge 0).
\end{displaymath}
From \eqref{30}, we have
\begin{displaymath}
\lim_{\lambda\rightarrow 0}W_{m,\lambda}(n,k)=W_{m}(n,k),
\end{displaymath}
where $W_{m}(n,k)$ are the Whitney numbers of the second kind given by
\begin{equation*}
(mx+1)^{n}=\sum_{k=0}^{n}W_{m}(n,k)m^{k}(x)_{k},\quad (n\ge 0),\quad (\mathrm{see}\ [16,17,28]).
\end{equation*}
We consider the degenerate Dowling polynomials defined by
\begin{equation}
D_{m,\lambda}(n,x)=\sum_{k=0}^{n}W_{m,\lambda}(n,k)x^{k},\quad (n\ge 0). \label{31}	
\end{equation}
Note that
\begin{displaymath}
D_{m}(n,x)=\lim_{\lambda\rightarrow 0}D_{m,\lambda}(n,x)=\sum_{k=0}^{n}W_{m}(n,k)x^{k}
\end{displaymath}
are the Dowling polynomials. \par
For $r\ge 0$, the degenerate $r$-Whitney numbers of the second kind are given by
\begin{equation}
(mx+r)_{n,\lambda}=\sum_{k=0}^{n}W_{m,\lambda}^{(r)}(n,k)m^{k}(x)_{k},\quad (n\ge 0),\quad (\mathrm{see}\ [16]).\label{32}
\end{equation}
In view of \eqref{31}, the degenerate $r$-Dowling polynomials are defined by
\begin{equation}
D_{m,\lambda}^{(r)}(n,x)=\sum_{k=0}^{n}W_{m,\lambda}^{(r)}(n,k)x^{k},\quad (n\ge 0),\quad (\mathrm{see}\ [16,17]).\label{33}
\end{equation}
From \eqref{15}, \eqref{30} and \eqref{32}, we note that
\begin{equation}
\big(m\fa+1\big)_{n,\lambda}=\sum_{k=0}^{n}W_{m,\lambda}(n,k)m^{k}\big(\ap\big)^{k}a^{k}, \label{34}
\end{equation}
and
\begin{equation}
\big(m\fa+r\big)_{n,\lambda}=\sum_{k=0}^{n}W_{m,\lambda}^{(r)}(n,k)m^{k}\big(\ap\big)^{k}a^{k},\quad (\mathrm{see}\ [16]).\label{35}	
\end{equation}
By \eqref{34}, we get
\begin{align}
\big(m\fa+1\big)_{l+n,\lambda}&=\big(m\fa+1-l\lambda\big)_{n,\lambda}\big(m\fa+1\big)_{l,\lambda}\label{36}\\
&=\big(m\fa+1-l\lambda\big)_{n,\lambda}\sum_{j=0}^{l}W_{m,\lambda}(l,j)m^{j}\big(\ap\big)^{j}a^{j}\nonumber\\
&=\sum_{j=0}^{l}W_{m,\lambda}(l,j)m^{j}\big(m\fa+1-l\lambda\big)_{n,\lambda}\big(\ap\big)^{j}a^{j},\nonumber
\end{align}
where $l,n$ are nonnegative integers. \par
By \eqref{36} and Theorem 2.3, we get
\begin{equation}
\big(m\fa+1\big)_{l+n,\lambda}=\sum_{j=0}^{l}W_{m,\lambda}(l,j)m^{j}\big(\ap\big)^{j}\big(m\fa+1+mj-l\lambda\big)_{n,\lambda}a^{j},\label{37}	
\end{equation}
where $l,n$ are nonnegative integers. \\
By Theorem 2.2, \eqref{31} and \eqref{34}, we get
\begin{align}
&\big(m\fa+1\big)_{l,\lambda}e^{\ap}\big|0\big\rangle=\sum_{j=0}^{l}W_{m,\lambda}(l,j)m^{j}\big(\ap\big)^{j}a^{j}e^{\ap}\big|0\big\rangle\label{38} \\
&=\sum_{j=0}^{l}W_{m,\lambda}(l,j)m^{j}\big(\ap\big)^{j}e^{\ap}\big|0\rangle=D_{m,\lambda}(l,m\ap)e^{\ap}\big|0\rangle, \nonumber
\end{align}
where $l$ is a nonnegative integer. \par
For $l,n\ge 0$, by Theorem 2.2, \eqref{26}, \eqref{37} and \eqref{38}, we get
\begin{align}
&D_{m,\lambda}(l+n,m\ap)e^{\ap}\big|0\big\rangle=\big(m\fa+1\big)_{l+n,\lambda}e^{\ap}\big|0\big\rangle \label{39}	\\
&=\sum_{j=0}^{l}W_{m,\lambda}(l,j)m^{j}\big(\ap\big)^{j}\big(m\fa+1+mj-l\lambda\big)_{n,\lambda}a^{j}e^{\ap}\big|0\big\rangle \nonumber\\
&=\sum_{j=0}^{l}W_{m,\lambda}(l,j)m^{j}\big(\ap\big)^{j}\big(m\fa+1+mj-l\lambda\big)_{n,\lambda}e^{\ap}\big|0\big\rangle \nonumber\\
&=\sum_{j=0}^{l}\sum_{k=0}^{n}W_{m,\lambda}(l,j)\binom{n}{k}m^{j}(\ap)^{j}(mj-l\lambda)_{n-k,\lambda}\big(m\fa+1\big)_{k,\lambda}e^{\ap}\big|0\big\rangle \nonumber\\
&=\sum_{j=0}^{l}\sum_{k=0}^{n}W_{m,\lambda}(l,j)\binom{n}{k}m^{j}(\ap)^{j}(mj-l\lambda)_{n-k,\lambda}D_{m,\lambda}(k,m\ap)e^{\ap}\big|0\big\rangle. \nonumber
\end{align}
Therefore, by \eqref{39}, we obtain the following theorem.
\begin{theorem}
For $n,l\ge 0$, we have
\begin{displaymath}
D_{m,\lambda}(n+l,m\ap)=\sum_{j=0}^{l}\sum_{k=0}^{n}\binom{n}{k}W_{m,\lambda}(l,j)m^{j}(\ap)^{j}(mj-l\lambda)_{n-k,\lambda}D_{m,\lambda}(k,m\ap),
\end{displaymath}
and hence
\begin{displaymath}
D_{m,\lambda}(n+l,x)=\sum_{j=0}^{l}\sum_{k=0}^{n}\binom{n}{k}W_{m,\lambda}(l,j) x^{j}(mj-l\lambda)_{n-k,\lambda}D_{m,\lambda}(k,x).
\end{displaymath}
\end{theorem}
From Theorem 2.2, \eqref{33} and \eqref{35}, we note that
\begin{align}
\big(m\fa+r\big)_{n,\lambda}e^{\ap}|0\rangle &=\sum_{k=0}^{n}W_{m,\lambda}^{(r)}(n,k)m^{k}\big(\ap\big)^{k}a^{k}e^{\ap}\big|0\big\rangle \label{40} \\
&=\sum_{k=0}^{n}W_{m,\lambda}^{(r)}(n,k)m^{k}\big(\ap\big)^{k}e^{\ap}\big|0\big\rangle\nonumber\\
&=D_{m,\lambda}^{(r)}\big(n,m\ap\big)e^{\ap}\big|0\big\rangle.\nonumber	
\end{align}
By Theorem 2.3, \eqref{26} and \eqref{35}, we get
\begin{align}
\big(m\fa+r\big)_{l+n,\lambda}&=\big(m\fa+r-l\lambda\big)_{n,\lambda}\big(m\fa+r\big)_{l,\lambda}\label{41}\\
&=\big(m\fa+r-l\lambda\big)_{n,\lambda}\sum_{j=0}^{l}W_{m,\lambda}^{(r)}(l,j)m^{j}\big(\ap\big)^{j}a^{j} \nonumber\\
&=\sum_{j=0}^{l}W_{m,\lambda}^{(r)}(l,j)m^{j}\big(m\fa+r-l\lambda\big
)_{n,\lambda}\big(\ap)^{j}a^{j} \nonumber\\
&=\sum_{j=0}^{l}W_{m,\lambda}^{(r)}(l,j)m^{j}\big(\ap\big)^{j}\big(m\fa+r+mj-l\lambda\big)_{n,\lambda}a^{j}\nonumber\\
&=\sum_{j=0}^{l}\sum_{k=0}^{n}W_{m,\lambda}^{(r)}(l,j)\binom{n}{k}m^{j}\big(\ap\big)^{j}(mj-l\lambda)_{n-k,\lambda}\big(m\fa+r\big)_{k,\lambda}a^{j}.\nonumber
\end{align}
By Theorem 2.2, \eqref{40} and \eqref{41}, we have
\begin{align*}
&D_{m,\lambda}^{(r)}(l+n,m\ap)e^{\ap}\big|0\big\rangle=\big(m\fa+r\big)_{l+n,\lambda}e^{\ap}\big|0\big\rangle \\
&=\sum_{j=0}^{l}\sum_{k=0}^{n}W_{m,\lambda}^{(r)}(l,j)\binom{n}{k}m^{j}\big(\ap\big)^{j}(mj-l\lambda)_{n-k,\lambda}\big(m\fa+r\big)_{k,\lambda}a^{j}e^{\ap}\big|0\big\rangle \\
&=\sum_{j=0}^{l}\sum_{k=0}^{n}W_{m,\lambda}^{(r)}(l,j)\binom{n}{k}m^{j}\big(\ap\big)^{j}(mj-l\lambda)_{n-k,\lambda}\big(m\fa+r\big)_{k,\lambda}e^{\ap}\big|0\big\rangle \\
	&=\sum_{j=0}^{l}\sum_{k=0}^{n}W_{m,\lambda}^{(r)}(l,j)\binom{n}{k}m^{j}\big(\ap\big)^{j}(mj-l\lambda)_{n-k,\lambda}D_{m,\lambda}^{(r)}(k,m\ap)e^{\ap}\big|0\big\rangle.
\end{align*}
This shows the following theorem.
\begin{theorem}
For $l,n\ge 0$, we have
\begin{equation*}
D_{m,\lambda}^{(r)}(l+n,m\ap)
= \sum_{j=0}^{l}\sum_{k=0}^{n}\binom{n}{k}W_{m,\lambda}^{(r)}(l,j)m^{j}\big(\ap\big)^{j}(mj-l\lambda)_{n-k,\lambda}D_{m,\lambda}^{(r)}(k,m\ap),
\end{equation*}
and hence
\begin{equation}
D_{m,\lambda}^{(r)}(l+n,x)=	\sum_{j=0}^{l}\sum_{k=0}^{n}\binom{n}{k}W_{m,\lambda}^{(r)}(l,j)x^{j}(mj-l\lambda)_{n-k,\lambda}D_{m,\lambda}^{(r)}(k,x). \label{42}
\end{equation}
\end{theorem}
From \eqref{42} and letting $\lambda \rightarrow 0$, we have
\begin{equation*}
D_{m}^{(r)}(l+n,x)=\sum_{j=0}^{l}\sum_{k=0}^{n}\binom{n}{k}W_{m}^{(r)}(l,j)x^{j}(mj)^{n-k}D_{m}^{(r)}(k,x),
\end{equation*}
where $D_{m}^{(r)}(k,x)$ are the $r$-Dowling polynomials given by
\begin{equation*}
D_{m}^{(r)}(n,x)=\sum_{k=0}^{n}W_{m}^{(r)}(n,k)x^{k},\quad (n\ge 0).
\end{equation*}
\section{Conclusion}
In recent years, various degenerate versions of many special polynomials and numbers have been explored with renewed interests. Studying degenerate versions was initiated by Carlitz in his work on degenerate Bernoulli and Euler polynomials (see [7]).
In this paper, we showed Spivey's type recurrence relations for the degenerate Bell polynomials, the degenerate Dowling polynomials and the degenerate $r$-Dowling polynomials, respectively by using normal orderings of degenerate $n$-th powers of $\ap a$, $m \ap a +1$, and $m \ap a +r$. \par
Indeed, we derived the following recurrence relations for the degenerate Dowling polynomials $D_{m,\lambda}(n,x)$ and the degenerate $r$-Dowling polynomials $D_{m,\lambda}^{(r)}(n,x)$:
\begin{align*}
&D_{m,\lambda}(n+l,x)=\sum_{j=0}^{l}\sum_{k=0}^{n}\binom{n}{k}W_{m,\lambda}(l,j) x^{j}(mj-l\lambda)_{n-k,\lambda}D_{m,\lambda}(k,x), \\
&D_{m,\lambda}^{(r)}(l+n,x)=	\sum_{j=0}^{l}\sum_{k=0}^{n}\binom{n}{k}W_{m,\lambda}^{(r)}(l,j)x^{j}(mj-l\lambda)_{n-k,\lambda}D_{m,\lambda}^{(r)}(k,x),
\end{align*}
respectively from
\begin{align*}
&\big(m\fa+1\big)_{n,\lambda}=\sum_{k=0}^{n}W_{m,\lambda}(n,k)m^{k}\big(\ap\big)^{k}a^{k}, \\
&\big(m\fa+r\big)_{n,\lambda}=\sum_{k=0}^{n}W_{m,\lambda}^{(r)}(n,k)m^{k}\big(\ap\big)^{k}a^{k},	
\end{align*}
where $W_{m,\lambda}(n,k)$ and $W_{m,\lambda}^{(r)}(n,k)$ are the degenerate Whitney numbers of the second kind and the degenerate $r$-Whitney numbers of the second kind. \par

\end{document}